\documentclass[a4paper]{article}

\usepackage{amsmath,amssymb,graphicx}
\usepackage{subfigure}

\usepackage{lscape}
\usepackage{mathrsfs}

\usepackage{algorithmic}
\usepackage[ruled]{algorithm}

\usepackage{multirow}
\usepackage{verbatim}
\usepackage{enumerate}
\usepackage{color}

\newcommand{\qed}{$\hfill{\Box}$}

\usepackage{longtable}

\newtheorem{remark}{Remark}

\title{ Singular Values of Riemann Curvature Tensor }

\author{
Xiaokai He$^a$\thanks{ Email: hexiaokai77@163.com.  }
\quad
Hua Xiang$^b$\thanks{
Corresponding author.
E-mail: hxiang@whu.edu.cn.   }
\\ \\
{\small $^a$ School of Mathematics and Computational Science,} \\ {\small Hunan First Normal University, Changsha 410205, China}\\
{\small $^b$ School of Mathematics and Statistics, Wuhan University, Wuhan 430072, China}
}

\begin{document}

\date{}
\maketitle

\begin{abstract}
We introduce the concept of singular values for the Riemann curvature tensor,  a central mathematical tool in Einstein's theory of general relativity. We study the properties related to the singular values, and investigate five typical cases to show its relationship to the Ricci scalar and other invariants.
\end{abstract}


\section{Introduction}

Let $M$ be the Riemannian manifold equipped with Riemannian metric $g$.
The curvature tensor for Levi-Civita connection $\nabla$ is called  Riemann curvature tensor, or Riemann tensor.
The (0,4) type Riemannian curvature tensor is a quadrilinear
mapping\cite{DoCarmo_book92}:
$$ R: ~ T_p M \times T_p M \times T_p M \times T_p M  \rightarrow \mathbb{R} , $$
$$ R(W, Z, X, Y) := \langle W, R(X,Y) Z  \rangle, \quad  \forall ~W, X, Y, Z \in T_p M .$$
where $R(X,Y) =  [\nabla_X, \nabla_Y] - \nabla_{[X,Y]}$, and $T_p M$ is the tangent space of $M$ at the point $p$.

Let $\{ \partial_i \}$ be the set of basis vector of a local coordinate, and $g_{ij} := g( \partial_i, \partial_j )$.
Using $[\partial_i, \partial_j] = 0$ and $\nabla_{\partial_i} \partial_k = \Gamma^l_{ik} \partial_l$, where $\Gamma^l_{ik}$ is the Christoffel symbols of Levi-Civita connection,  we can calculate that
\begin{eqnarray*}
R( \partial_i, \partial_j ) \partial_k = ( \partial_i \Gamma^l_{jk} - \partial_j \Gamma^l_{ik}   + \Gamma^h_{jk}\Gamma^l_{ih} - \Gamma^h_{ik} \Gamma^l_{jh}) \partial_l  \equiv  R^l_{kij} \partial_l  ,
\end{eqnarray*}
$$ R(\partial_i, \partial_j, \partial_k, \partial_l) = \langle  \partial_i,  R(\partial_k,\partial_l) \partial_j \rangle = g_{ih} R^h_{jkl} \equiv  R_{ijkl}. $$

Let $W = w^i \partial_i$, $X = x^j \partial_j$, $Y = y^k \partial_k$, $Z = z^l \partial_l$.
We can verify that
$$R(Y,Z) X  = x^j y^k z^l  R( \partial_k, \partial_l)  \partial_j =  x^j y^k z^l R^h_{jkl} \partial_h, $$
$$R(W,X,Y,Z) = \langle W, R(Y,Z) X \rangle  = w^i x^j y^k z^l R_{ijkl}. $$

The curvature tensor has the following symmetry properties:
$$R_{ijkl} = -R_{jikl} = -R_{ijlk} = R_{klij}, $$
and the first (or algebraic) Bianchi identity:
$$   R_{ijkl} + R_{iljk} + R_{iklj} =0 . $$
For the 4D case, there are only twenty independent entries among 256 components due to these symmetries \cite{BaezMuniain_Book1994,Dirac_Book75,MTW1973}.

The contraction yields the Ricci tensor $R_{ik}$ and the Ricci scalar $R$ as follows.
$$ R_{ik}  = R^m_{imk} = g^{hj} R_{hijk} , \qquad
 R =  R^k_k = g^{ik} R_{ik} . $$

Recently, the tensor eigenvalues problem has been studied
extensively\cite{Qi_JSC05,QiChenChen18}. In \cite{XiangQiWei_CMS18}
we investigate the M-eigenvalue of Riemann curvature tensor. In this
paper, we will introduce and discuss the singular value of Riemann
curvature tensor. For the singular values of matrices and integral
operators, one can refer to the review \cite{Stewart_SIREV93}, where
Stewart survey the contributions of five mathematicians who
established and developed the singular value decomposition: E.
Beltrami, C. Jordan, J. J. Sylvester, E. Schmidt, H. Weyl. The first
three came to it through linear algebra, while the last two
approached it from integral equations.


In the following of the paper, we introduce the definition of singular value problem for Riemann tensor and investigate its properties in section 2. Then in section 3, we study five typical cases, calculate the singular values and examine their relationship to the well-known invariants.

\section{Singular value problem of Riemann tensor}

We recall the singular value problem (SVP) of a real matrix $A$ by starting with
 a bilinear form $f(x,y) = x^T A y$, and seek the maximum and minimum of $f$ subject to $||x|| = ||y|| =1$,
resulting in the following equations:
\begin{eqnarray*}
A y     &=& \sigma     x, \\
x^T A   &=& \lambda    y^T.
\end{eqnarray*}
It is obvious that $\lambda = \sigma$.

We introduce the singular value problem of Riemann curvature tensor by considering the following optimization problem on Riemann manifold:
\begin{eqnarray*}
&& \min_{W,X,Y,Z \in T_p M }  R(W,X,Y,Z)            \\
&& \text{s.t.} \quad \langle W,W \rangle = \langle X,X \rangle = \langle Y,Y \rangle = \langle Z,Z \rangle = 1.
\end{eqnarray*}
The feasible set of this optimization problem is compact.  Hence, it always has global optimal solutions.
 It has only equality constraints.   By optimization theory, its optimal Lagrangian multipliers  $(\theta,\lambda,\mu,\nu)$   always exist and are real.
\begin{eqnarray*}
 L(W,X,Y,Z, \sigma, \lambda, \mu, \nu) &=& R_{ijkl} w^i x^j y^k z^l  + \theta (g_{ij} w^i w^j -1) \\
 && + \lambda (g_{ij} x^i x^j -1) + \mu (g_{ij} y^i y^j -1) + \nu (g_{ij} z^i z^j -1) .
\end{eqnarray*}
The optimality condition reads
\begin{eqnarray*}
R_{ijkl} x^j y^k z^l = - 2 \theta     w_i,  &&\quad
R_{ijkl} w^i y^k z^l = - 2 \lambda    x_j,  \\
R_{ijkl} w^i x^j z^l = - 2 \mu        y_k,  &&\quad
R_{ijkl} w^i x^j y^k = - 2 \nu        z_l.
\end{eqnarray*}
This ensures the existence of such $\theta$, $\lambda$, $\mu$ and $\nu$, and they are real.
It is easy to verify that $\theta = \lambda = \mu =\nu $.

We may rewrite the singular value problem  of Riemann tensor as seeking the normalized singular vectors $W, X, Y, Z$ and the singular value $\sigma$ in a coordinate-free manner.

In an abstract form, the SVP of Riemann tensor is defined by the 5-tuple $(W,X,Y,Z, \sigma)$ satisfying
the following formulae:
\begin{eqnarray}
 R(Y,Z) X  &=&    \sigma W, \label{s1}\\
 R(Z,Y) W  &=&    \sigma X, \label{s2}\\
 R(W,X) Z  &=&    \sigma Y, \label{s3}\\
 R(X,W) Y  &=&    \sigma Z, \label{s4}
 \end{eqnarray}
with \begin{eqnarray}\langle W,W \rangle= \langle X,X \rangle=
\langle Y,Y \rangle= \langle Z,Z \rangle=1.      \label{s5}
\end{eqnarray}

The componentwise SVP of \eqref{s1}--\eqref{s4} reads:
\begin{eqnarray*}
R^i_{jkl} x^j y^k z^l &=& \sigma w^i,  \\
R^i_{jkl} w^j z^k y^l &=& \sigma x^i,  \\
R^i_{jkl} z^j w^k x^l &=& \sigma y^i,  \\
R^i_{jkl} y^j x^k w^l &=& \sigma z^i.
\end{eqnarray*}

The SVP can be introduced in another way. Consider the following associated Lagrangian function \cite{LHLim_IEEE05}.
$$ \mathcal{L}: \quad T_p M \times T_p M \times T_p M \times T_p M \times \mathbb{R} \longrightarrow \mathbb{R}. $$
$$ \mathcal{L} ( W, X, Y, Z, \sigma ) := R (W, X, Y, Z) -  \sigma  \left(
\langle W,W \rangle^\frac{1}{2} \langle X,X \rangle^\frac{1}{2}
\langle Y,Y \rangle^\frac{1}{2} \langle Z,Z \rangle^\frac{1}{2}  - 1 \right) . $$

In component form it reads
$$  \mathcal{L}  = R_{ijkl} w^i x^j y^k z^l -  \sigma  \left(
(g_{ij} w^i w^j)^\frac{1}{2}  (g_{kl} x^k x^l)^\frac{1}{2} (g_{mn}
y^m y^n)^\frac{1}{2}  (g_{pq} z^p z^q)^\frac{1}{2}  - 1 \right) ,$$
from which the SVP of Riemann tensor can be derived. For example,
$\partial_{w^i} \mathcal{L} = 0 $ yields
$$  R_{ijkl}   x^j y^k z^l -  \sigma   (
g_{ij} w^i w^j)^{-\frac{1}{2}} (g_{ij} w^j)  (g_{kl} x^k x^l)^\frac{1}{2}
 (g_{mn} y^m y^n)^\frac{1}{2}  (g_{pq} z^p z^q)^\frac{1}{2} = 0.  $$
That is,
$$   R_{ijkl}   x^j y^k z^l  ( g_{kl} x^k x^l)^{-\frac{1}{2}}  (g_{mn} y^m y^n)^{-\frac{1}{2}}
(g_{pq} z^p z^q)^{-\frac{1}{2}} =   \sigma  w_i   (g_{ij} w^i w^j)^{-\frac{1}{2}}   .  $$
With a bit abuse of the notation, we use $W$ to denote the normalized vector $\langle W,W \rangle^{-\tfrac{1}{2}} W$.
Then $\partial_W L = 0$ gives $R_{ijkl} x^j y^k z^l = \sigma w_i$, that is, $R(Y,Z) X  = \sigma W$.

The partial derivatives with respect to $X, Y, Z$ give the similar results. We summarize them in the following.
\begin{eqnarray*}
 \partial_W L = 0 ~\Longrightarrow \eqref{s1}, && \quad  \partial_X L = 0 ~\Longrightarrow \eqref{s2}, \\
 \partial_Y L = 0 ~\Longrightarrow \eqref{s3}, && \quad  \partial_Z L = 0 ~\Longrightarrow \eqref{s4}.
\end{eqnarray*}

\begin{remark}
We introduce the SVP from the optimization problem over the Riemannian manifold, where the feasible set of the optimization problem is compact. The definition should not be limited to a Riemannian manifold, and we can extend the definition to a pseudo-Riemannian manifold. But for the pseudo-Riemannian manifold, there may be no maximum
or minimum of the optimization problem, since the constraint set is not compact any more.
But the introduction via the Lagrangian function $\mathcal{L}$ does not need such extra consideration.
\end{remark}

\begin{remark}
On the Riemann manifold, we can check that $\sigma = R(W,X,Y,Z)$. For the Lorentz manifold, if $\langle W,W \rangle = \langle X,X \rangle =\langle Y,Y \rangle = 1$, and $\langle Z,Z \rangle = -1$, we have $\sigma = 0$, since the inner product with the formulae \eqref{s1}--\eqref{s3} yields $ R(W,X,Y,Z)=\sigma $ while the formula \eqref{s4} gives $ R(W,X,Y,Z)= -\sigma $.
\end{remark}

\begin{remark}
We can check that the following 5-tuples
$$ ( X, X, Y, Z, \sigma=0 ) , $$
$$ ( W, X, Z, Z, \sigma=0 ) , $$
$$ ( W, W, W, W, \sigma=0 )   $$
solve the singular value problem trivially, and hence $\sigma=0$ always be a singular value for
the problem.
\end{remark}

\begin{remark}
 The formulae \eqref{s1}--\eqref{s4} can be rewritten as :
$$ R(Y,Z) (W+iX) = i \sigma (W+iX), \quad R(W,X) (Y+iZ) = i \sigma (Y+iZ). $$
We can clearly see the symmetry between the pairs $(W,X)$ and $(Y,Z)$.

When $W = Y$, $X=Z$, we have the M-eigenvalue problem
\cite{XiangQiWei_CMS18,HanDaiQi09,QDH09,XiangQiWei2017}
$$ R(Y,Z) Z = \sigma Y, \quad  R(Z,Y) Y = \sigma Z. $$
That is, $ R(Y,Z) (Y+iZ) = i \sigma (Y+iZ) $.
\end{remark}

%

One can obtain the following properties of this singular value
problem.

\textbf{Proposition 1.} If $W,X,Y,Z$ and $\sigma$ solve the singular
value problem $(\ref{s1})\sim (\ref{s5})$, and $\sigma$ is nonzero, then
$$ \langle W,X \rangle = 0,  \qquad  \langle Y,Z  \rangle = 0. $$

\textbf{Proposition 2.} If $ (W,X,Y,Z)$ and $\sigma$ solve the SVP $(\ref{s1})\sim (\ref{s5})$, then the following 5-tuples also solve the SVP $(\ref{s1})\sim (\ref{s5})$.
$$ ( -W,X,Y,Z, -\sigma ) ,$$
$$ ( W,-X,Y,Z, -\sigma ) ,$$
$$ ( W,X,-Y,Z, -\sigma ) ,$$
$$ ( W,X,Y,-Z, -\sigma ) .$$

\textbf{Proposition 3.} If $ (W,X,Y,Z)$ and $\sigma$ solve the SVP $(\ref{s1})\sim (\ref{s5})$, then the following 5-tuples also solve the SVP.
$$ ( X,W,Y,Z, -\sigma ) ,$$
$$ ( W,X,Z,Y, -\sigma ) ,$$
$$ ( X,W,Z,Y,  \sigma ) ,$$
%
and
$$ ( Y,Z,W,X, \sigma ) .$$

{\bf Corollary 1.}
From the Propostion 2, we find that
$$ ( -W,-X,Y,Z, \sigma  ),$$
$$ ( W,X,-Y,-Z, \sigma  ), $$
$$ ( W,-X,-Y,Z, \sigma ), $$
$$ ( W,-X,Y,-Z, \sigma ), $$
$$ ( -W,X,-Y,Z, \sigma ), $$
$$ ( -W,X,Y,-Z, \sigma ), $$
$$ ( W,-X,-Y,-Z, -\sigma ), $$
$$ ( -W,X,-Y,-Z, -\sigma ), $$
$$ ( -W,-X,Y,-Z, -\sigma ), $$
$$ ( -W,-X,-Y,Z, -\sigma ), $$
and
$$ ( -W,-X,-Y,-Z, \sigma )  $$
 also solve the singular value problem $(\ref{s1})\sim(\ref{s5})$.

Applying the Propostion 3, we have that
$$ ( Z,Y,W,X, -\sigma ), $$
$$ ( Y,Z,X,W, -\sigma ), $$
$$ ( Z,Y,X,W,  \sigma )  $$
 also solve the singular value problem $(\ref{s1})\sim(\ref{s5})$.

\textbf{Proposition 4.} If $W,X,Y,Z$ and $\sigma$ solve the singular
value problem $(\ref{s1})\sim (\ref{s5})$, then
$$ \left( \tfrac{1}{\sqrt{2}}(W-X),\ \tfrac{1}{\sqrt{2}}(W+X),\ Y,Z,
~\sigma \right) , $$
$$ \left( W,X, \tfrac{1}{\sqrt{2}}(Y-Z),\ \tfrac{1}{\sqrt{2}}(Y+Z),\
~\sigma \right) , $$ and
$$ \left(  \tfrac{1}{\sqrt{2}}(W+X),\ \tfrac{1}{\sqrt{2}}(W-X),\ \tfrac{1}{\sqrt{2}}(Y+Z),\
\tfrac{1}{\sqrt{2}}(Y-Z),\  \sigma \right)   $$
also solve the singular value problem $(\ref{s1})\sim(\ref{s5})$.

We can write out more 5-tuples that satisfying the SVP $(\ref{s1})\sim(\ref{s5})$   like those in  Proposition 2, Proposition 2 and Corollary 1. For simplicity, we omit them here.

\begin{remark}
Proposition 2 can be derived from the definition of SVP directly, without using the symmetries of Riemann tensor. From Proposition 2, we find that the singular value $\sigma$ can be restricted to be always non-negative.
\end{remark}



\section{Case study}

\textbf{Example 1.} Singular value problem of the 2 dimensional
standard sphere.

\begin{eqnarray}
ds^2=d\theta^2+\sin^2\theta d\varphi^2
\end{eqnarray}

The nonzero components of Riemann tensor read

$$R_{1212}=-R_{1221}=R_{2121}=-R_{2112}=\sin^2\theta . $$

Suppose that
\begin{eqnarray}
( W = w^i \partial_i, X = x^j \partial_j, Y = y^k
\partial_k, Z = z^l \partial_l, \ \sigma )
\end{eqnarray}
solve the singular value problem $(\ref{s1})\sim(\ref{s5})$.

Then we have
\begin{eqnarray*}
\sigma w^1=\sin^2\theta  x^2(y^1z^2-y^2z^1),\\
\sigma w^2=x^1(y^2z^1-y^1z^2),\\
\sigma x^1=\sin^2\theta w^2(y^2z^1-y^1z^2),\\
\sigma x^2=w^1(y^1z^2-y^2z^1),\\
\sigma y^1=\sin^2\theta z^2(w^1x^2-w^2x^1),\\
\sigma y^2=z^1(w^2x^1-w^1x^2),\\
\sigma z^1=\sin^2\theta y^2(w^2x^1-w^1x^2),\\
\sigma z^2 =y^1(w^1x^2-w^2x^1).
\end{eqnarray*}

One can check that
\begin{eqnarray}
W = \pm \frac{1}{\sin\theta}\partial_{\varphi}, \
X = \pm \partial_{\theta},\
Y = \pm \frac{1}{\sin\theta}\partial_{\varphi}, \
Z = \pm \partial_{\theta},\
\sigma = 1
\end{eqnarray}
solves the singular value problem.

According to the properties in section 2, we can write down more,
for example,
$W = Y = \pm \partial_{\theta}$, $X=Z= \pm \frac{1}{\sin\theta}\partial_{\varphi}$. We omit them here for simplicity.

{\bf Example 2.}
In the following we investigate  the Riemannian manifold with constant sectional curvature $\kappa$.
$\forall ~ W, X, Y, Z \in T_pM $ we have \cite[P.149]{Lee_Book97}
$$ R(W,X,Y,Z) = \langle W, R(X,Y)Z \rangle = \kappa ( \langle W,X \rangle
 \langle Z,Y \rangle - \langle W,Y \rangle \langle Z,X \rangle ), $$
equivalently,
$$ R(X,Y) Z =  \kappa ( \langle Z,Y \rangle X - \langle Z,X \rangle Y ). $$

The corresponding SVP reads as follows.
$$ R(Y,Z) X = \kappa ( \langle Z,X \rangle Y - \langle Y,X \rangle Z ) =    \sigma W, $$
$$ R(Z,Y) W = \kappa ( \langle Y,W \rangle Z - \langle Z,W \rangle Y ) =    \sigma X, $$
$$ R(W,X) Z = \kappa ( \langle X,Z \rangle W - \langle W,Z \rangle X ) =    \sigma Y, $$
$$ R(X,W) Y = \kappa ( \langle W,Y \rangle X - \langle X,Y \rangle W ) =    \sigma Z. $$
We can derive that
\begin{equation} \label{eqn:deriveEqnSigma}
\kappa ( \langle Z,X \rangle \langle Y,W \rangle - \langle Y,X \rangle \langle Z,W \rangle ) =    \sigma .
\end{equation}

For the case of $\kappa=0$, one can obtain that $\sigma=0$
identically, and $\forall ~ W, X, Y, Z \in T_pM $ solve the singular
value problem.

For the case of $\kappa\neq 0$, we know that $\sigma=0$ is a
singular value of the singular value problem from the remark 3. Besides, by proposition 1, when $\sigma\neq 0$,
$ \langle W,X  \rangle = 0 $,
$ \langle Y,Z  \rangle = 0 $.
Using these, we have
\begin{eqnarray}
  \kappa \langle Z,X  \rangle &=& \sigma   \langle  W,Y \rangle,  \label{eqn:deriveFourEqns1} \\
 -\kappa \langle Y,X  \rangle &=& \sigma   \langle  W,Z \rangle, \label{eqn:deriveFourEqns2}\\
 -\kappa \langle Z,W  \rangle &=& \sigma   \langle  X,Y \rangle,  \label{eqn:deriveFourEqns3} \\
  \kappa \langle Y,W  \rangle &=& \sigma   \langle  X,Z \rangle. \label{eqn:deriveFourEqns4}
\end{eqnarray}

From these four identities \eqref{eqn:deriveFourEqns1}--\eqref{eqn:deriveFourEqns4}, we obtain
$$ \kappa \sigma \langle  W,Y \rangle^2 = \kappa \sigma \langle  X,Z \rangle^2 ,
\qquad  - \kappa \sigma \langle  W,Z \rangle^2 = - \kappa \sigma \langle  X,Y \rangle^2 . $$
That is,
$$  | \langle  W,Y \rangle | = | \langle  X,Z \rangle | ,
\qquad  | \langle  W,Z \rangle | = | \langle  X,Y \rangle | . $$

Using \eqref{eqn:deriveFourEqns1}--\eqref{eqn:deriveFourEqns4}, we can also obtain
$$  ( \kappa^2 - \sigma^2 ) \langle  W,Y \rangle \langle  X,Z \rangle = 0,
\qquad  ( \kappa^2 - \sigma^2 )  \langle  W,Z \rangle  \langle  X,Y \rangle = 0. $$
Taking into account the results above, we have
$$  ( \kappa^2 - \sigma^2 ) \langle  W,Y \rangle^2   = 0,
\qquad  ( \kappa^2 - \sigma^2 )  \langle  W,Z \rangle^2    = 0. $$

Substituting \eqref{eqn:deriveFourEqns1}--\eqref{eqn:deriveFourEqns2} into \eqref{eqn:deriveEqnSigma}, we have
$$ \langle  W,Y \rangle^2 + \langle  W,Z \rangle^2 = 1 . $$
And hence,
$\langle  W,Y \rangle$ and $\langle  W,Z \rangle$ cannot be zero simutaneously.
So we have $\kappa^2 - \sigma^2 = 0$. That is,
$$ \sigma =    | \kappa |  . $$

{\bf Example 3.} It's commonly known that  the Riemann tensor can be splitted into trace and tracefree parts.
\begin{eqnarray*}
R_{ijkl} =
C_{ijkl} - \tfrac{1}{n-2} \left( R_{il} g_{jk} - R_{ik} g_{jl} + g_{il} R_{jk} - g_{ik} R_{jl} \right)
\\  + \tfrac{R}{(n-1)(n-2)} \left( g_{il} g_{jk} - g_{ik} g_{jl} \right) ,
\end{eqnarray*}
where $C_{ijkl}$ is the tracefree Weyl tensor. For the conformally flat manifold, the Weyl tensor $C_{ijkl}=0$, then the singular value problem reads
\begin{eqnarray}
- \tfrac{1}{n-2} \left( \langle x,y \rangle R_{il} z^l   - \langle x,z \rangle R_{ik} y^k +  ( R_{jk} x^j y^k ) z_i - ( R_{jl} x^j z^l ) y_i   \right)     \nonumber \\
+ \tfrac{R}{(n-1)(n-2)} \left( \langle x,y \rangle z_i  - \langle x,z \rangle  y_i \right)
  = \sigma w_i,     \label{eqn:Conformal_w} \\
- \tfrac{1}{n-2} \left( \langle w,z \rangle R_{jk} y^k   - \langle w,y \rangle R_{jl} z^l +  ( R_{il} w^i z^l ) y_j - ( R_{ik} w^i y^k  ) z_j   \right)        \nonumber \\
+ \tfrac{R}{(n-1)(n-2)} \left( \langle w,z \rangle y_j  - \langle w,y \rangle  z_j \right)
  = \sigma x_j,     \label{eqn:Conformal_x} \\
- \tfrac{1}{n-2} \left( \langle w,z \rangle R_{jk} x^j   - \langle x,z \rangle R_{ik} w^i +  ( R_{il} w^i z^l ) x_k - ( R_{jl} x^j z^l  ) w_k   \right)        \nonumber \\
+ \tfrac{R}{(n-1)(n-2)} \left( \langle w,z \rangle x_k  - \langle x,z \rangle  w_k \right)
  = \sigma y_k,     \label{eqn:Conformal_y}  \\
- \tfrac{1}{n-2} \left( \langle x,y \rangle R_{il} w^i   - \langle w,y \rangle R_{jl} x^j +  ( R_{jk} x^j y^k ) w_l - ( R_{ik} w^i y^k ) x_l   \right)         \nonumber \\
+ \tfrac{R}{(n-1)(n-2)} \left( \langle x,y \rangle w_l  - \langle w,y \rangle  x_l \right)
  = \sigma z_l .     \label{eqn:Conformal_z}
\end{eqnarray}

We seek the nonzero singular value. We can easily check that
$ \langle w,x \rangle  = \langle y,z \rangle  = 0 $. Taking the inner product and using the unit norm condition, we have
\begin{eqnarray}
- \tfrac{1}{n-2} \left(  (R_{il} z^l w^i) \langle x,y \rangle  - (R_{ik}  y^k w^i)  \langle x,z \rangle +  ( R_{jk} y^k x^j  ) \langle w,z \rangle - ( R_{jl} z^l x^j  ) \langle w,y \rangle  \right)         \nonumber
\\
+ \tfrac{R}{(n-1)(n-2)} \left( \langle w,z \rangle  \langle x,y \rangle    - \langle w,y \rangle  \langle x,z \rangle  \right) = \sigma  . \label{eqn:sigmaExpression0}
\end{eqnarray}
In the following, we would like to derive some explicit expressions related to $\sigma$.
We will derive the expressions for the terms $R_{jl} x^j z^l$, $R_{jk} x^j y^k$, $R_{il} w^i z^l$ and $R_{ik} w^i y^k$ respectively.

The inner product of \eqref{eqn:Conformal_w} with $y$ yields
\begin{eqnarray*}
- \tfrac{1}{n-2} \left( \langle x,y \rangle R_{il} z^l y^i  - \langle x,z \rangle R_{ik} y^k y^i   - R_{jl} x^j z^l 
\right)
  - \tfrac{R}{(n-1)(n-2)}   \langle x,z \rangle  
  = \sigma \langle w,y \rangle,
\end{eqnarray*}
and the inner product of \eqref{eqn:Conformal_y} with $w$ gives
\begin{eqnarray*}
- \tfrac{1}{n-2} \left( - \langle x,z \rangle R_{ik} w^i w^k  + \langle w,z \rangle R_{jk} x^j w^k   - R_{jl} x^j z^l  
\right)
  - \tfrac{R}{(n-1)(n-2)}   \langle x,z \rangle  
  = \sigma \langle y,w \rangle .
\end{eqnarray*}
From the two equations above we have
\begin{equation} \label{eqn:Conformal_derived1}
 ( R_{ik} w^i w^k - R_{ik} y^k y^i ) \langle x,z \rangle = R_{jk}x^j w^k \langle w,z \rangle - R_{il} z^l y^i \langle y,x \rangle ,
\end{equation}
$$
\tfrac{1}{n-2} R_{jl} x^j z^l   = \sigma \langle w,y \rangle  + \tfrac{R}{(n-1)(n-2)}   \langle x,z \rangle +  \tfrac{1}{n-2}  \left( \langle x,y \rangle R_{il} z^l y^i  - \langle x,z \rangle R_{ik} y^k y^i  \right) .
$$

Using the inner product of \eqref{eqn:Conformal_w} with $z$, and \eqref{eqn:Conformal_z} with $w$, we have
\begin{eqnarray*}
- \tfrac{1}{n-2} \left( \langle x,y \rangle R_{il} z^l z^i   - \langle x,z \rangle R_{ik} y^k z^i +   R_{jk} x^j y^k \right)     + \tfrac{R}{(n-1)(n-2)}  \langle x,y \rangle    = \sigma \langle w,z \rangle  ,
\end{eqnarray*}
$$
- \tfrac{1}{n-2} \left( \langle x,y \rangle R_{il} w^i w^l   - \langle w,y \rangle R_{jl} x^j w^l +    R_{jk} x^j y^k  \right)  + \tfrac{R}{(n-1)(n-2)}  \langle x,y \rangle    = \sigma \langle z,w \rangle.
$$
Then we obtain
\begin{equation} \label{eqn:Conformal_derived2}
 ( R_{il} w^i w^l - R_{il} z^l z^i ) \langle x,y \rangle = R_{jl}x^j w^l \langle w,y \rangle - R_{ik} y^k z^i  \langle z,x \rangle ,
\end{equation}
$$
- \tfrac{1}{n-2} R_{jk} x^j y^k   = \sigma \langle w,z \rangle -  \tfrac{R}{(n-1)(n-2)}  \langle x,y \rangle  + \tfrac{1}{n-2}   \left( \langle x,y \rangle R_{il} z^l z^i   - \langle x,z \rangle R_{ik} y^k z^i     \right)    .
$$

Using the inner product of \eqref{eqn:Conformal_x} with $y$, and \eqref{eqn:Conformal_y} with $x$, we have
\begin{eqnarray*}
- \tfrac{1}{n-2} \left( \langle w,z \rangle R_{jk} y^k y^j  - \langle w,y \rangle R_{jl} z^l y^j +   R_{il} w^i z^l   \right)  + \tfrac{R}{(n-1)(n-2)}  \langle w,z \rangle    = \sigma \langle x,y \rangle,
\end{eqnarray*}
$$
- \tfrac{1}{n-2} \left( \langle w,z \rangle R_{jk} x^j x^k   - \langle x,z \rangle R_{ik} w^i x^k +  R_{il} w^i z^l  \right)  + \tfrac{R}{(n-1)(n-2)} \langle w,z \rangle   = \sigma \langle y,x \rangle.
$$
Hence,
\begin{equation} \label{eqn:Conformal_derived3}
 ( R_{jk} x^j x^k - R_{jk} y^k y^j ) \langle w,z \rangle =  R_{ik} w^i x^k  \langle x,z \rangle  -  R_{jl} z^l y^j \langle y,w \rangle ,
\end{equation}
$$
- \tfrac{1}{n-2} R_{il} w^i z^l   = \sigma \langle x,y \rangle - \tfrac{R}{(n-1)(n-2)}  \langle w,z \rangle  +  \tfrac{1}{n-2}   \left( \langle w,z \rangle R_{jk} y^k y^j  - \langle w,y \rangle R_{jl} z^l y^j   \right).
$$

Similarly, using the inner product of \eqref{eqn:Conformal_x} with $z$, and \eqref{eqn:Conformal_z} with $x$, we have
\begin{eqnarray*}
- \tfrac{1}{n-2} \left( \langle w,z \rangle R_{jk} y^k z^j  - \langle
w,y \rangle R_{jl} z^l z^j -  R_{ik} w^i y^k  \right)     -
\tfrac{R}{(n-1)(n-2)}  \langle w,y \rangle  = \sigma \langle x,z
\rangle,
\end{eqnarray*}
$$
- \tfrac{1}{n-2} \left( \langle x,y \rangle R_{il} w^i x^l   - \langle w,y \rangle R_{jl} x^j x^l  -  R_{ik} w^i y^k \right) - \tfrac{R}{(n-1)(n-2)}  \langle w,y \rangle    = \sigma \langle z,x \rangle .
$$

Therefore,
\begin{equation} \label{eqn:Conformal_derived4}
 ( R_{jl} x^j x^l - R_{jl} z^l z^j ) \langle w,y \rangle =  R_{il} w^i x^l  \langle x,y \rangle  -  R_{jk} y^k z^j \langle z,w \rangle ,
\end{equation}
$$
\tfrac{1}{n-2} R_{ik} w^i y^k   = \sigma \langle x,z \rangle + \tfrac{R}{(n-1)(n-2)}  \langle w,y \rangle  + \tfrac{1}{n-2}  \left( \langle w,z \rangle R_{jk} y^k z^j  - \langle w,y \rangle R_{jl} z^l z^j \right)  .
$$


Substituting the above expressions for $R_{il} w^i z^l$, $R_{ik} w^i y^k$,  $R_{jk} x^j y^k$ and $R_{jl} x^j z^l$  into \eqref{eqn:sigmaExpression0}, we have
\begin{eqnarray}
\sigma  ( \langle w,y \rangle^2 + \langle w,z \rangle^2 + \langle x,y \rangle^2 + \langle x,z \rangle^2  - 1 ) =  \nonumber \\
 ( \langle w,y \rangle \langle x,z \rangle - \langle x,y \rangle \langle w,z \rangle  )
 \tfrac{1}{n-2}  \left( R_{jk} y^j y^k + R_{jl} z^l z^j  - \tfrac{R}{n-1} \right). \label{eqn:sigmaExpression1}
\end{eqnarray}
Similarly, we can derive that
\begin{eqnarray}
\sigma  ( \langle w,y \rangle^2 + \langle w,z \rangle^2 + \langle x,y \rangle^2 + \langle x,z \rangle^2 - 1 ) =  \nonumber \\
   ( \langle w,y \rangle \langle x,z \rangle - \langle x,y \rangle \langle w,z \rangle )
 \tfrac{1}{n-2} \left( R_{ik} w^i w^k + R_{jl} x^j x^l  - \tfrac{R}{n-1} \right). \label{eqn:sigmaExpression2}
\end{eqnarray}
A by-product of \eqref{eqn:sigmaExpression1} and \eqref{eqn:sigmaExpression2} is that  $R_{ij} w^i w^j + R_{ij} x^i x^j = R_{ij} y^i y^j + R_{ij} z^i z^j$.
It's hard to further simplify and derive an explicit expression for $\sigma$. In the following we consider three special cases to simplify the expression.

 \textbf{Case 1}. For the M-eigenvalue problem where $W=Y$ and $X=Z$, we can check that  \eqref{eqn:sigmaExpression1} and \eqref{eqn:sigmaExpression2} reduce to
$$
\sigma  =
  \tfrac{1}{n-2}
 \left( R_{ik} w^i w^k + R_{jl} x^j x^l  - \tfrac{R}{n-1} \right).
$$

  \textbf{Case 2}. Suppose that $(\lambda, y)$ and $(\mu, z)$
are the eigen-pairs of the Ricci tensor. That is, $R_{ij} y^j =
\lambda y_i$ and $R_{ij} z^j = \mu z_i$. Let $w = \pm y$ and $x =
\pm z$. Then we have
$$\sigma =  \tfrac{1}{n-2} \left(\lambda + \mu - \tfrac{R}{n-1} \right) . $$

 \textbf{Case 3}. For Einstein manifold, we have
\begin{eqnarray*}
R_{ij}=\kappa g_{ij} .
\end{eqnarray*}
Then \eqref{eqn:Conformal_w} yields
\begin{eqnarray*}
\tfrac{1}{n-2}\left(2\kappa- \tfrac{R}{n-1}\right)\langle x,z \rangle
 =\sigma \langle w,y \rangle,
\end{eqnarray*}
and \eqref{eqn:Conformal_x} yields
\begin{eqnarray*}
\tfrac{1}{n-2}\left(2\kappa-\tfrac{R}{n-1}\right)\langle w,y \rangle
 =\sigma \langle x,z \rangle .
\end{eqnarray*}
And hence one can deduce that
\begin{eqnarray*}
 \sigma= \tfrac{1}{n-2}\left(2\kappa -\tfrac{R}{n-1}\right) .
\end{eqnarray*}


{\bf Example 4. } We consider the  Schwarzschild solution.
In the Schwarzschild coordinates, with signature $(-1, 1, 1, 1)$, the line element for the Schwarzschild metric has the form
$$ds^2 = - f(r)  dt^2 + f(r)^{-1} dr^2 + r^2 ( d\theta^2 + \sin^2 \theta d\phi^2 ) ,$$
where  $f(r) = 1 - \frac{r_s}{r}$ and $r_s = {2 M} $. Here we adopt the geometric unit $c=G=1.$

Define  $A=\frac{ M f(r)}{r^3} =\frac{ r_s }{2 r^3}  f(r)$, $B = \frac{ M}{ r^3 f(r) } = \frac{ r_s }{ 2 r^3 f(r) }  $, $C=\frac{M }{ r} =\frac{r_s }{2 r}$, $D=\frac{M \sin^2 \theta}{ r} =\frac{r_s }{2 r} \sin^2 \theta  $. The nonzero components are expressed by
\begin{eqnarray*}
 R^0_{101} =-R^0_{110}= 2B, \qquad R^0_{220}=-R^0_{202}= C, \qquad R^0_{330}=-R^0_{303} = D, \\
 R^1_{001}=-R^1_{010} = 2A, \qquad R^1_{221}=-R^1_{212}= C, \qquad R^1_{331}=-R^1_{313} = D ,\\
 R^2_{020} =-R^2_{002}= A,  \qquad R^2_{112}=-R^2_{121}= B, \qquad R^2_{323} =-R^2_{332}= 2D,\\
 R^3_{030}=-R^3_{003} = A,  \qquad R^3_{113}=-R^3_{131}= B, \qquad R^3_{232}=-R^3_{223} = 2C .
\end{eqnarray*}

We have the following system of polynomial equations with seventeen variables $(w^0, w^1, w^2, w^3, x^0, x^1, x^2, x^3,  y^0, y^1, y^2, y^3, z^0, z^1, z^2, z^3, \sigma)$.
\begin{eqnarray*}
&2 B x^1 (y^0z^1-y^1z^0)  +C  x^2 (y^2z^0-y^0 z^2)  +D    x^3 (y^3z^0-y^0 z^3) & = \sigma w^0, \\
&2 A x^0 (y^0z^1-y^1z^0)  +C  x^2 (y^2z^1-y^1 z^2)  +D    x^3 (y^3z^1-y^1 z^3) & = \sigma w^1, \\
&  A x^0 (y^2z^0-y^0z^2)  +B  x^1 (y^1 z^2-y^2z^1)  +2D   x^3 (y^2 z^3-y^3z^2) & = \sigma w^2, \\
&  A x^0 (y^3z^0-y^0z^3)  +B  x^1 (y^1 z^3-y^3z^1)  +2C   x^2 (y^3z^2-y^2 z^3) & = \sigma w^3, \\
& 2B w^1 (z^0y^1-z^1y^0)  +Cw^2 (z^2y^0-z^0 y^2)  +D w^3 (z^3y^0-z^0 y^3) & = \sigma x^0, \\
& 2A w^0 (z^0y^1-z^1y^0)  +Cw^2 (z^2y^1-z^1 y^2)  +D w^3 (z^3y^1-z^1 y^3) & = \sigma x^1, \\
&  A w^0 (z^2y^0-z^0y^2)  +Bw^1 (z^1 y^2-z^2y^1) +2D w^3 (z^2 y^3-z^3y^2) & = \sigma x^2, \\
&  A w^0 (z^3y^0-z^0y^3)  +Bw^1 (z^1 y^3-z^3y^1) +2C w^2 (z^3y^2-z^2 y^3) & = \sigma x^3, \\
& 2B z^1 (w^0x^1-w^1x^0)  +Cz^2 (w^2x^0-w^0 x^2)  +D z^3 (w^3x^0-w^0 x^3) & = \sigma y^0, \\
& 2A z^0 (w^0x^1-w^1x^0)  +Cz^2 (w^2x^1-w^1 x^2)  +D z^3 (w^3x^1-w^1 x^3) & = \sigma y^1, \\
&  A z^0 (w^2x^0-w^0x^2)  +Bz^1 (w^1 x^2-w^2x^1) +2D z^3 (w^2 x^3-w^3x^2) & = \sigma y^2, \\
&  A z^0 (w^3x^0-w^0x^3)  +Bz^1 (w^1 x^3-w^3x^1) +2C z^2 (w^3x^2-w^2 x^3) & = \sigma y^3, \\
& 2B y^1 (x^0w^1-x^1w^0)  +Cy^2 (x^2w^0-x^0 w^2)  +D y^3 (x^3w^0-x^0 w^3) & = \sigma z^0, \\
& 2A y^0 (x^0w^1-x^1w^0)  +Cy^2 (x^2w^1-x^1 w^2)  +D y^3 (x^3w^1-x^1 w^3) & = \sigma z^1, \\
&  A y^0 (x^2w^0-x^0w^2)  +By^1 (x^1 w^2-x^2w^1) +2D y^3 (x^2 w^3-x^3w^2) & = \sigma z^2, \\
&  A y^0 (x^3w^0-x^0w^3)  +By^1 (x^1 w^3-x^3w^1) +2C y^2 (x^3w^2-x^2 w^3) & = \sigma z^3,
\end{eqnarray*}
together with the four constraints:
\begin{eqnarray*}
&    g_{00} w^0 w^0 + g_{11} w^1 w^1 + g_{22} w^2 w^2 + g_{33} w^3 w^3 & =  \pm 1 , \\
&    g_{00} x^0 x^0 + g_{11} x^1 x^1 + g_{22} x^2 x^2 + g_{33} x^3 x^3 & =  \pm 1 , \\
&    g_{00} y^0 y^0 + g_{11} y^1 y^1 + g_{22} y^2 y^2 + g_{33} y^3 y^3 & =  \pm 1 , \\
&    g_{00} z^0 z^0 + g_{11} z^1 z^1 + g_{22} z^2 z^2 + g_{33} z^3 z^3 & =  \pm 1 .
\end{eqnarray*}

Let $S^{ij} = y^i z^j - y^j z^i$.
From the first four equalities, we have
$$
\begin{bmatrix}
0 & 2S^{01} & S^{20} & S^{30} \\ 2S^{01} & 0 & S^{21} & S^{31} \\
S^{20} & S^{12} & 0 & 2S^{23} \\ S^{30} & S^{13} & 2S^{32} & 0 \\
\end{bmatrix}
\begin{bmatrix}
A & 0 & 0 & 0 \\ 0 & B & 0 & 0 \\ 0 & 0 & C & 0 \\ 0 & 0 & 0 & D \\
\end{bmatrix}
\begin{bmatrix}
x^0 \\ x^1 \\ x^2 \\ x^3
\end{bmatrix} = \sigma
\begin{bmatrix}
w^0 \\ w^1 \\ w^2 \\ w^3
\end{bmatrix}
$$
For simplicity, we denote it by
$$ S  \Lambda x = \sigma w , $$
where $\Lambda :=  \text{diag}(A,B,C,D)  =  \frac{r_s} {2 r^3} \text{diag}(c^2 f(r), f(r)^{-1}, r^2, r^2 \sin^2 \theta) $, and $S$ can be defined correspondingly.

Similarly, from the next four equalities, we have
$$ S \Lambda w = - \sigma x .$$


Let $T^{ij} = w^i x^j - w^j x^i$.
From the third four equalities, we have
$$
\begin{bmatrix}
0 & 2T^{01} & T^{20} & T^{30} \\ 2T^{01} & 0 & T^{21} & T^{31} \\
T^{20} & T^{12} & 0 & 2T^{23} \\ T^{30} & T^{13} & 2T^{32} & 0 \\
\end{bmatrix}
\begin{bmatrix}
A & 0 & 0 & 0 \\ 0 & B & 0 & 0 \\ 0 & 0 & C & 0 \\ 0 & 0 & 0 & D \\
\end{bmatrix}
\begin{bmatrix}
z^0 \\ z^1 \\ z^2 \\ z^3
\end{bmatrix} = \sigma
\begin{bmatrix}
y^0 \\ y^1 \\ y^2 \\ y^3
\end{bmatrix}
$$
For simplicity, we denote it by
$$ T  \Lambda z = \sigma y . $$
Similarly, from the last four equalities, we have
$$ T \Lambda y = - \sigma z .$$


Rewriting these four equations in a compact form, we have
$$
\begin{bmatrix}
 0 & S \Lambda & 0 & 0 \\ - S \Lambda & 0 & 0 & 0 \\
 0 & 0 & 0 & T \Lambda \\ 0 & 0 & - T \Lambda & 0
 \end{bmatrix}
\begin{pmatrix} w \\ x \\ y \\ z \end{pmatrix} = \sigma
\begin{pmatrix} w \\ x \\ y \\ z \end{pmatrix} .
$$
Note that
$$
\begin{pmatrix} \sigma & - S \Lambda \\   S \Lambda & \sigma \end{pmatrix}  =
\begin{pmatrix} I & 0 \\   \sigma^{-1} S \Lambda & I \end{pmatrix}
\begin{pmatrix} \sigma & - S \Lambda \\   0 & \sigma + \sigma^{-1} (S \Lambda)^2 \end{pmatrix}, $$
$$ \det
\begin{pmatrix} \sigma & - S \Lambda \\   S \Lambda & \sigma \end{pmatrix} =
\det \left(  \sigma^2   + (S \Lambda)^2 \right) = \det \left( S \Lambda + i \sigma I   \right) \det \left( S \Lambda - i \sigma I   \right) ,  $$

$$ \det
\begin{pmatrix}
 \sigma & -S \Lambda & 0 & 0 \\  S \Lambda & \sigma & 0 & 0 \\
 0 & 0 & \sigma & -T \Lambda \\ 0 & 0 & T \Lambda & \sigma
 \end{pmatrix}
= \det \left(  \sigma^2   + (S \Lambda)^2 \right) \det \left(  \sigma^2  + (T \Lambda)^2 \right)
. $$
The nonzero solution requires that
$$\det \left( S \Lambda \pm i \sigma I   \right) =0, \quad \text{or} \quad \det \left( T \Lambda \pm i \sigma I   \right) =0. $$

 Direct calculation yields that
\begin{eqnarray*}
\det \left( S  \pm i \sigma \Lambda^{-1}   \right) =   \det(S) +
[ \sigma^2 + A B (2S^{01})^2 + A C (S^{20})^2 + A D (S^{30})^2  \\
 - C D (2S^{23})^2 - B D (S^{13})^2 - B C (S^{21})^2   ] \sigma^2  \det(\Lambda^{-1})  .
\end{eqnarray*}
A sufficient condition for the existence of nonzero solution is that the two terms in the right hand size of above formula are zeros.

For the first term, we can verify that
$$ \det(S) = - ( 2 S^{01} \cdot 2 S^{23} + S^{20} S^{13} + S^{30} S^{21} )^2 .$$
Hence, $ \det(S) = 0$ is equivalent to 
$4 S^{01} S^{23} + S^{20} S^{13} + S^{30} S^{21} =0 $.
Substituting the express for $S^{ij}$, we have $y^0 z^1 = y^1 z^0$ or $y^2 z^3=y^3 z^2$, that is, $S^{01} S^{23} =0$.
Similarly, we have $T^{01} T^{23} =0$ that is, $w^0 x^1 = w^1 x^0$ or $w^2 x^3 = w^3 x^2$.

We choose
\begin{eqnarray*}
 y^0 = z^0 , \quad y^1 = z^1, && y^2 = -z^2, \quad y^3 = -z^3,  \\
 x^0 = w^0 , \quad x^1 = w^1, && x^2 = -w^2, \quad x^3 = -w^3,
\end{eqnarray*}
such that $y \neq z$ and $w \neq x$.
Under such settings, the SVP reduces to the following polynomial system involving the unknowns $x$, $y$ and $\sigma$.

\begin{eqnarray*}
& C  x^2 y^0 y^2  +D    x^3 y^0 y^3 & = \tfrac{\sigma}{2} x^0, \\
& C  x^2 y^1 y^2  +D    x^3 y^1 y^3 & = \tfrac{\sigma}{2} x^1, \\
&-A  x^0 y^0 y^2  +B    x^1 y^1 y^2 & = \tfrac{\sigma}{2} x^2, \\
&-A  x^0 y^0 y^3  +B    x^1 y^1 y^3 & = \tfrac{\sigma}{2} x^3, \\
& C  y^2 x^0 x^2  +D    y^3 x^0 x^3 & = \tfrac{\sigma}{2} y^0, \\
& C  y^2 x^1 x^2  +D    y^3 x^1 x^3 & = \tfrac{\sigma}{2} y^1, \\
&-A  y^0 x^0 x^2  +B    y^1 x^1 x^2 & = \tfrac{\sigma}{2} y^2, \\
&-A  y^0 x^0 x^3  +B    y^1 x^1 x^3 & = \tfrac{\sigma}{2} y^3,
\end{eqnarray*}
together with the four constraints:
\begin{eqnarray*}
\langle x, x \rangle & = & \tfrac{2r^3}{r_s} [- A(x^0)^2 + B(x^1)^2 + C(x^2)^2 + D(x^3)^2] =   1, \\
\langle y, y \rangle & = & \tfrac{2r^3}{r_s} [- A(y^0)^2 + B(y^1)^2 + C(y^2)^2 + D(y^3)^2] =   1.
\end{eqnarray*}
Solving it, we have
$$ \sigma =  \frac{r_s}{2 r^3}  = \pm \frac{M}{  r^3} = \pm \sqrt{ \frac{K_1}{48} } , $$
where
$K_1 = R_{abcd} R^{abcd} = 12 \left( \frac{r_s}{r^3} \right)^2 = \frac{48  M^2 }{  r^6} $ is the Kretschman curvature invariant.
Note that the constraints $\langle x, x \rangle = \langle y, y \rangle = -1$ give the same results, while $\langle x, x \rangle =1$ and $ \langle y, y \rangle = -1$ yield $\sigma=0$.

\qed

\textbf{Example 5.}
In the Boyer-Lindquist coordinate, the Kerr solution reads
\begin{eqnarray}
ds^2=-\bigg{(}1-\frac{2Mr}{r^2+a^2\cos^2\theta}\bigg{)}\bigg{(}dt+
\frac{2Mar\sin^2\theta}{r^2-2Mr+a^2\cos^2\theta}d\varphi\bigg{)}^2\nonumber\\
+\bigg{(}1-\frac{2Mr}{r^2+a^2\cos^2\theta}\bigg{)}^{-1}\bigg{[}(r^2-2Mr+a^2)\sin^2\theta d\varphi^2\nonumber\\
+(r^2-2Mr+a^2\cos^2\theta)\bigg{(}\frac{1}{r^2-2Mr+a^2}dr^2+d\theta^2\bigg{)}\bigg{]} .
\end{eqnarray}
By choosing the following Newman-Penrose null tetrads
\cite{Chandrasekhar_Book83}
\begin{eqnarray*}
l^i = \bigg{(}\frac{r^2+a^2}{r^2-2Mr+a^2},1,0,\frac{a}{r^2-2Mr+a^2}\bigg{)},\\
n^i = \bigg{(}\frac{r^2+a^2}{2(r^2+a^2\cos^2\theta)},-\frac{r^2-2Mr+a^2}{2(r^2+a^2\cos^2\theta)},
0,\frac{a}{2(r^2+a^2\cos^2\theta)}\bigg{)},\\
m^i = \frac{1}{\sqrt{2}(r+ia\cos\theta)}\bigg{(}ia\sin\theta,0,1,\frac{i}{\sin\theta}\bigg{)},\\
\bar{m}^i = \frac{1}{\sqrt{2}(r-ia\cos\theta)}\bigg{(}-ia\sin\theta,0,1,\frac{-i}{\sin\theta}\bigg{)},
\end{eqnarray*}
one can obtain that
\begin{eqnarray}
\Psi_0=\Psi_1=\Psi_3=\Psi_4=0,\ \ \Psi_2 = \frac{M}{(r-ia\cos\theta)^3}.
\end{eqnarray}
Let
\begin{eqnarray*}
W &=& w^1 l+w^2n+w^3m+w^4\bar{m},\\
X &=& x^1 l+x^2n+x^3m+x^4\bar{m},\\
Y &=& y^1 l+y^2n+y^3m+y^4\bar{m},\\
Z &=& z^1 l+z^2n+z^3m+z^4\bar{m},
\end{eqnarray*}
then the singular value problem becomes a system of polynomial equations with seventeen variables $(w^0, w^1, w^2, w^3, x^0, x^1, x^2, x^3,  y^0, y^1, y^2, y^3, z^0, z^1, z^2, z^3, \sigma)$.
\begin{eqnarray*}
 -(\Psi_2 + \bar{\Psi}_2) (x^2 y^1 z^2 - x^2 y^2 z^1) + (\Psi_2 - \bar{\Psi}_2) (x^2 y^3 z^4 - x^2 y^4 z^3)   \\
    \qquad + \Psi_2 (x^3 y^2 z^4 - x^3 y^4 z^2) + \bar{\Psi}_2 (x^4 y^2 z^3 - x^4 y^3 z^2) = - \sigma w^2,             \\
 (\Psi_2 + \bar{\Psi}_2) (x^1 y^1 z^2 - x^1 y^2 z^1) - (\Psi_2 - \bar{\Psi}_2) (x^1 y^3 z^4 - x^1 y^4 z^3)    \\
   \qquad + \bar{\Psi}_2 (x^3 y^1 z^4 - x^3 y^4 z^1) + \Psi_2 (x^4 y^1 z^3 - x^4 y^3 z^1) = - \sigma w^1,             \\
   (\Psi_2 - \bar{\Psi}_2) (x^4 y^1 z^2 - x^4 y^2 z^1) - (\Psi_2 + \bar{\Psi}_2) (x^4 y^3 z^4 - x^4 y^4 z^3)  \\
   \qquad  -\Psi_2 (x^1 y^2 z^4 - x^1 y^4 z^2) - \bar{\Psi}_2 (x^2 y^1 z^4 - x^2 y^4 z^1) = \sigma w^4,             \\
 - (\Psi_2 - \bar{\Psi}_2) (x^3 y^1 z^2 - x^3 y^2 z^1) + (\Psi_2 + \bar{\Psi}_2) (x^3 y^3 z^4 - x^3 y^4 z^3) \\
  \qquad   -\bar{\Psi}_2 (x^1 y^2 z^3 - x^1 y^3 z^2) - \Psi_2 (x^2 y^1 z^3 - x^2 y^3 z^1) = \sigma w^3,             \\
  \end{eqnarray*}
  \begin{eqnarray*}
 -(\Psi_2+\bar{\Psi}_2)(w^2z^1y^2-w^2z^2y^1)+(\Psi_2-\bar{\Psi}_2)(w^2z^3y^4-w^2z^4y^3) \\
  \qquad +\Psi_2(w^3z^2y^4-w^3z^4y^2)+\bar{\Psi}_2(w^4z^2y^3-w^4z^3y^2)=-\sigma x^2, \\
 (\Psi_2+\bar{\Psi}_2)(w^1z^1y^2-w^1z^2y^1)-(\Psi_2-\bar{\Psi}_2)(w^1z^3y^4-w^1z^4y^3)\\
 \qquad +\bar{\Psi}_2(w^3z^1y^4-w^3z^4y^1)+\Psi_2(w^4z^1y^3-w^4z^3y^1)=-\sigma x^1,\\
   (\Psi_2-\bar{\Psi}_2)(w^4z^1y^2-w^4z^2y^1)-(\Psi_2+\bar{\Psi}_2)(w^4z^3y^4-w^4z^4y^3) \\
 \qquad -\Psi_2(w^1z^2y^4-w^1z^4y^2)-\bar{\Psi}_2(w^2z^1y^4-w^2z^4y^1) =\sigma x^4,\\
  -(\Psi_2-\bar{\Psi}_2)(w^3z^1y^2-w^3z^2y^1)+(\Psi_2+\bar{\Psi}_2)(w^3z^3y^4-w^3z^4y^3) \\
 \qquad -\bar{\Psi}_2(w^1z^2y^3-w^1z^3y^2)-\Psi_2(w^2z^1y^3-w^2z^3y^1) =\sigma x^3,\\
\end{eqnarray*}
\begin{eqnarray*}
-(\Psi_2+\bar{\Psi}_2)(z^2w^1x^2-z^2w^2x^1)+(\Psi_2-\bar{\Psi}_2)(z^2w^3x^4-z^2w^4x^3)\\
    +\Psi_2(z^3w^2x^4-z^3w^4x^2)+\bar{\Psi}_2(z^4w^2x^3-z^4w^3x^2)      =-\sigma y^2,\\
(\Psi_2+\bar{\Psi}_2)(z^1w^1x^2-z^1w^2x^1)-(\Psi_2-\bar{\Psi}_2)(z^1w^3x^4-z^1w^4x^3)\\
    +\bar{\Psi}_2(z^3w^1x^4-z^3w^4x^1)+{\Psi}_2(z^4w^1x^3-z^4w^3x^1)    =-\sigma y^1,\\
(\Psi_2-\bar{\Psi}_2)(z^4w^1x^2-z^4w^2x^1)-(\Psi_2+\bar{\Psi}_2)(z^4w^3x^4-z^4w^4x^3)\\
    -\Psi_2(z^1w^2x^4-z^1w^4x^2)-\bar{\Psi}_2(z^2w^1x^4-z^2w^4x^1)      =\sigma y^4,\\
-(\Psi_2-\bar{\Psi}_2)(z^3w^1x^2-z^3w^2x^1)+(\Psi_2+\bar{\Psi}_2)(z^3w^3x^4-z^3w^4x^3)\\
    -\bar{\Psi}_2(z^1w^2x^3-z^1w^3x^2)-{\Psi}_2(z^2w^1x^3-z^2w^3x^1)    =\sigma y^3,\\
\end{eqnarray*}
\begin{eqnarray*}
-(\Psi_2+\bar{\Psi}_2)(y^2x^1w^2-y^2x^2w^1)+(\Psi_2-\bar{\Psi}_2)(y^2x^3w^4-y^2x^4w^3)\\
    +\Psi_2(y^3x^2w^4-y^3x^4w^2)+\bar{\Psi}_2(y^4x^2w^3-y^4x^3w^2)=-\sigma z^2, \\
(\Psi_2+\bar{\Psi}_2)(y^1x^1w^2-y^1x^2w^1)-(\Psi_2-\bar{\Psi}_2)(y^1x^3w^4-y^1x^4w^3)\\
    +\bar{\Psi}_2(y^3x^1w^4-y^3x^4w^1)+{\Psi}_2(y^4x^1w^3-y^4x^3w^1)=-\sigma z^1,\\
(\Psi_2-\bar{\Psi}_2)(y^4x^1w^2-y^4x^2w^1)-(\Psi_2+\bar{\Psi}_2)(y^4x^3w^4-y^4x^4w^3) \\
    -\Psi_2(y^1x^2w^4-y^1x^4w^2)-\bar{\Psi}_2(y^2x^1w^4-y^2x^4w^1) =\sigma z^4,\\
-(\Psi_2-\bar{\Psi}_2)(y^3x^1w^2-y^3x^2w^1)+(\Psi_2+\bar{\Psi}_2)(y^3x^3w^4-y^3x^4w^3) \\
    -\bar{\Psi}_2(y^1x^2w^3-y^1x^3w^2)-{\Psi}_2(y^2x^1w^3-y^2x^3w^1) =\sigma z^3.\\
\end{eqnarray*}

It hard to solve the above system of polynomial equations generally, so we will seek a special solution of this system.
Setting
$$
Z=X, \quad W=-Y
$$
and
$$
y^1=x^1,\ y^2=x^2,\ y^3=-x^3,\ y^4=-x^4,
$$
then the system of equations reduces to
\begin{eqnarray*}
2x^3x^4(\Psi_2+\bar{\Psi}_2)-\sigma=0,\\
2x^1x^2(\Psi_2+\bar{\Psi}_2)+\sigma=0.
\end{eqnarray*}
Combining with
\begin{eqnarray*}
\langle x,x \rangle =-2x^1x^2+2x^3x^4=1,
\end{eqnarray*}
we can obtain the singular value as
\begin{eqnarray}
\sigma=\frac{(r^3-3r\cos\theta^2)M}{(r^2+a^2\cos\theta^2)^3}=\sqrt{\frac{ |I| +  \Re(I) }{6}},
\end{eqnarray}
where $I\equiv \Psi_0\Psi_4-4\Psi_1\Psi_3+3\Psi_2^2$ is an invariant \cite{Stephani_Book2003}.

For Schwarzschild case, $I=\frac{3M^2}{r^6}$, and the singular value reduces to
\begin{eqnarray}
\sigma= \frac{M}{r^3}.
\end{eqnarray}

\section*{Acknowledgement}
 H. Xiang is supported by the National
Natural Science Foundation of China under grants 11571265, 11471253 and NSFC-RGC 
No.11661161017.  X. He was supported by the Grant of NSF of Hunan
No.2018JJ2073 and NSFC No.11401199.

\end{document}